\documentclass[11pt]{article}                  
\usepackage{amsmath,amsthm,amsfonts,graphicx}
 
\newtheorem{theorem}{Theorem}

\theoremstyle{definition}   
 
\newtheorem{definition}[theorem] {Definition}

\newtheorem{remark}[theorem]{Remark}

\DeclareMathOperator{\rank}{rank}
\DeclareMathOperator{\nrank}{nrank}

\usepackage{multirow}   
\begin{document} 
\title{On the Superfast Multipole Method}
\author{Victor Y. Pan} 
\author{Victor Y. Pan$^{[1, 2],[a]}$ and
John Svadlenka$^{[2],[c]}$ 
\and\\
$^{[1]}$ Department of Computer Science \\
Lehman College of the City University of New York \\
Bronx, NY 10468 USA \\
$^{[2]}$ Ph.D. Programs in  Computer Science and Mathematics \\
The Graduate Center of of the City University of New York \\
New York, NY 10036 USA \\
$^{[a]}$ victor.pan@lehman.cuny.edu \\ 
http://comet.lehman.cuny.edu/vpan/  \\
$^{[b]}$  jsvadlenka@gradcenter.cuny.edu \\ 
}

\date{}  
\maketitle  
  
  
\begin{abstract} 
  

   We call matrix algorithms {\em superfast} if they involve much fewer memory cells and flops than the input matrix has entries. Using such algorithms is indispensable for Big Data Mining and Analysis,   
   where the input matrices    are so immense that realistically one can only access a small fraction of all their entries.  A natural remedy is Low Rank Approximation of these 
   matrices,\footnote{Hereafter we use the acronym LRA.} which is routinely computed by means of 
   Cross--Approximation\footnote{Hereafter we use the acronym C–A.} iterations for more than a decade of worldwide application in computational practice.
We point out and extensively test an
     important application of superfast LRA to significant acceleration of the celebrated Fast Multipole Method,
    which  turns it into Supefast Multipole Method.
        
   
\end{abstract}
 
 




\section{Introduction: Superfast LRA}\label{sintro}


{\em Low rank approximation} (hereafter
{\em LRA}) of a matrix  is a fundamental subject of Numerical Linear Algebra  
and Computer Science. An $m\times n$ matrix $M$ admits its close approximation\footnote{Here and hereafter the concepts ``low", ``large", ``small", ``far", ``close",  etc.
 are defined in context. The inequalities
  $a\ll b$ and $b\gg a$ show that $|a/b|$ is small in context.} 
 of rank $r$ if and only if
 the matrix $M$ has {\em numerical rank} $r$ (then we write $\nrank(W)=r$) or equivalently if
 and only if
\begin{equation}\label{eqlrk}  
M=AB+E,~||E||/||M||\le \epsilon,
\end{equation}
for a  small integer $r$, 
$A\in \mathbb C^{m\times r}$,
 $B\in \mathbb C^{r\times n}$, a~fixed matrix~norm $||\cdot||$, and~a~small tolerance $\epsilon$.
Such an LRA approximates the $mn$  entries of $M$ by using $(m+n)r$ entries of $A$ and $B$ instead of $mn$ entries of $M$, and one can operate with a low rank matrix, e.g., multiply it by a vector superfast. This is a crucial benefit 
in applications of LRA to Big 
Data Mining and Analysis, where the input matrices $M$, e.g., unfolding matrices of multidimensional tensors,
 are so immense that realistically one can only access a tiny fraction of all their entries.
  LRA is a natural remedy, and
 for more than a decade  the 
Cross--Approximation (C-A) 
   iterations have routinely  been computing accurate LRA superfast in worldwide computational practice
  ( cf. \cite{T96}, \cite{GTZ97},
\cite{T00}, \cite{B00},
\cite{BR03},  \cite{GOSTZ10},
 \cite{B11},  \cite{O18},  \cite{OZ18}).
   

\section{An Application -- Superfast Multipole Method}\label{sextpr} 
 

Superfast LRA algorithms can be extended to numerous important computational problems linked to LRA.
Next we we point out and
    extensively test a
    simple but apparently unnoticed application of superfast LRA to significant acceleration of the celebrated Fast Multipole Method (FMM),
    which  turns it into {\em Superfast Multipole Method}. 
 

\subsection{Fast and Superfast Multipole Method}

%

FMM enables superfast  multiplication by a vector of so
called  HSS matrices
provided that low rank generators are available for its off-diagonal blocks.
 Such generators are not available in some important applications, however (see, e.g, \cite{XXG12}, 
 \cite{XXCB14}, and \cite{P15}), but  C--A algorithms  compute them
superfast, thus turning FMM into Superfast 
 Multipole Method.  
 Since the method is highly important
we supply some details of its  bottleneck stage of HSS computations, which we perform superfast by incorporating superfast LRA.                                                                                                                                                                                                               
 
 HSS matrices
 naturally extend the class of banded matrices and their inverses, 
are closely linked to FMM,  
 and 
are increasingly popular
(see \cite{BGH03},  \cite{GH03}, \cite{MRT05},
 \cite{CGS07},
 \cite{VVGM05}, 
\cite{VVM07/08}, 
 \cite{B10}, \cite{X12},   \cite{XXG12},
\cite{EGH13}, \cite{X13},
 \cite{XXCB14},
  and the bibliography therein). 

\begin{definition}\label{defneut} {\rm (Neutered Block Columns. See 
\cite{MRT05}.)} With each diagonal block of a block matrix 
associate 
its complement in its block column,
and call this complement a {\em neutered block column}.
\end{definition}

\begin{definition}\label{defqs} {\rm (HSS matrices. See 
\cite{CGS07},
 \cite{X12},  \cite{X13}, \cite{XXCB14}.)}
  
A block 
matrix $M$ of size  $m\times n$ is 
called an $r$-{\em HSS matrix}, for a positive integer $r$, 

(i) if all diagonal blocks of this matrix
consist of $O((m+n)r)$ entries overall
and
 
(ii) if $r$ is the maximal rank of its neutered block columns. 
\end{definition}

\begin{remark}\label{reqs}
Many authors work with $(l,u)$-HSS
(rather than $r$-HSS) matrices $M$ for which $l$ and $u$
are the maximal ranks of the sub- and super-diagonal blocks,
respectively.
The $(l,u)$-HSS and $r$-HSS matrices are closely related. 
If a neutered block column $N$
is the union of a sub-diagonal block $B_-$ and 
a super-diagonal block $B_+$,
then
 $\rank (N)\le \rank (B_-)+\rank (B_+)$,
 and so
an $(l,u)$-HSS matrix is an $r$-HSS matrix,
for $r\le l+u$,
while clearly an $r$-HSS matrix is  
a $(r,r)$-HSS matrix.
\end{remark}

The FMM exploits the $r$-HSS structure of a matrix as follows:

(i) Cover all off-block-diagonal entries
with a set of  non-overlapping neutered block columns.  

(ii) Express every neutered block column $N$ of this set
  as the product
  $FH$ of two  {\em generator
matrices}, $F$ of size $h\times r$
and $H$ of size $r\times k$. Call the 
pair $\{F,H\}$ a {\em length $r$ generator} of the 
neutered block column $N$. 

(iii)  Multiply 
the matrix $M$ by a vector by separately multiplying generators
and diagonal blocks by subvectors,  involving $O((m+n)r)$ flops
overall, and

(iv) in a more advanced application of  
 FMM solve a nonsingular $r$-HSS linear system of $n$
equations  by using
$O(nr\log^2(n))$ flops under some mild additional assumptions on  the input. 

This approach is readily extended to the same operations with  
$(r,\xi)$-{\em HSS matrices},
that is, matrices approximated by $r$-HSS matrices
within a perturbation norm bound $\xi$ where a  positive tolerance 
$\xi$ is small in context (for example, is the unit round-off).
 Likewise, one defines 
an  $(r,\xi)$-{\em HSS representation} and 
$(r,\xi)$-{\em generators}.

$(r,\xi)$-HSS matrices (for $r$ small in context)
appear routinely in matrix computations,
and computations with such matrices are 
performed  efficiently by using the
above techniques.


In some applications of the FMM (see  \cite{BGP05}, \cite{VVVF10})
stage (ii) is omitted because short generators for all 
neutered block columns are readily available,
 but this is not the case in a variety of other  important applications
 (see \cite{XXG12}, \cite{XXCB14}, and \cite{P15}). 
This stage of the computation of  generators is precisely 
 LRA of the neutered block
columns, which turns out to be 
the bottleneck stage of FMM in these applications, and superfast LRA algorithms  
provide a remedy. 

Indeed apply a fast  algorithm at this
 stage, e.g., the algorithm of \cite{HMT11}
 with a Gaussian multiplier.
Multiplication of a $q\times h$ matrix
by an $h\times r$ Gaussian matrix requires $(2h-1)qr$ flops,
while
standard HSS-representation of an $n\times n$ 
HSS matrix includes $q\times h$ neutered 
 block columns for $q\approx m/2$ and $h\approx n/2$. In this case 
the cost of computing an $r$-HSS representation of the
matrix $M$ has at least order $mnr$.
For $r\ll \min\{m,n\}$, this
is much  greater than 
 $O((m+n)r\log^2(n))$ flops, used  at the other stages of 
the computations. 

We alleviate such a problem, however, 
when we
compute LRA 
of  
$(r,\xi)$-generators
by applying superfast algorithms.
    

\section{Computation of  LRAs for benchmark HSS matrices}
\label{sHSS}
 

In this section, the contribution of the secind author, we cover our   tests of the Superfast Multipole Method where we  applied C--A iterations
in order to compute  LRA of the generators of the off-diagonal blocks of 
HSS matrices.  Namely we tested HSS matrices that approximate $1024 \times 1024$ Cauchy-like matrices 
derived from  benchmark Toeplitz matrices B, C, D, E, and F of \cite[Section 5]{XXG12}. For the computation of  LRA we applied the algorithm of \cite{GOSTZ10}.

Table 1
displays the relative errors of
the approximation of the  $1024 \times 1024$ HSS input matrices
 in the spectral and Chebyshev norms averaged over 100 tests.
 Each approximation was obtained by 
 means of   
 combining the exact diagonal blocks 
 and  LRA of the off-diagonal blocks.
  We
 computed  LRA of all these blocks 
 superfast.
 
In good accordance with extnsive empirical evidence about the power of C--A iterations, already the first C--A loop have consistently yielded reasonably close  LRA, but our  further improvement was achieved in five C--A loops in our tests for all but one of the five families of input matrices.

The reported HSS rank is the larger of the numerical ranks for the $512 \times 512$ off-diagonal blocks. This HSS rank was used as an upper bound in our binary search that determined the numerical rank of each off-diagonal block for the purpose of computing its LRA. We based 
the binary search on  minimizing the difference (in the spectral norm) between each off-diagonal block and its LRA. 

The output error norms were quite low. Even in the case of
 the matrix C, obtained from Prolate Toeplitz matrices -- extremely ill-conditioned, they ranged from $10^{-3}$
 to $10^{-6}$.

We have also performed further numerical experiments on all the HSS input matrices by using a hybrid LRA algorithm: we used random pre-processing with Gaussian and Hadamard (abridged and permuted)  multipliers by incorporating
 Algorithm 4.1 of \cite{HMT11}, but  only for the  off-diagonal blocks of smaller sizes while retaining our previous 
 way for computing  LRA of the larger off-diagonal blocks. We have not displayed the results of these experiments because they yielded no substantial improvement in accuracy in comparison to the exclusive use of the less expensive  LRA on all off-diagonal blocks. 

\begin{table}[ht]
\begin{center}
\begin{tabular}{|c|c|c|c|c|c|c|}
\hline
 &  & & \multicolumn{2}{c|}{\bf Spectral Norm} & \multicolumn{2}{c|}{\bf Chebyshev Norm} \\ \hline

{\bf Inputs} & {\bf C--A loops}  &  {\bf HSS rank} & {\bf mean} & {\bf std} & {\bf mean} & {\bf std} \\ \hline


\multirow{2}*{B}

\multirow{2}*{B}
&1 &  26 & 8.11e-07 & 1.45e-06 & 3.19e-07 & 5.23e-07 \\  \cline{2-7}
&5 &  26 & 4.60e-08 & 6.43e-08 & 7.33e-09 & 1.22e-08 \\  \hline
\multirow{2}*{C}
&1 &  16 & 5.62e-03 & 8.99e-03 & 3.00e-03 & 4.37e-03 \\  \cline{2-7}
&5 &  16 & 3.37e-05 & 1.78e-05 & 8.77e-06 & 1.01e-05 \\  \hline
\multirow{2}*{D}
&1 &  13 & 1.12e-07 & 8.99e-08 & 1.35e-07 & 1.47e-07 \\  \cline{2-7}
&5 &  13 & 1.50e-07 & 1.82e-07 & 2.09e-07 & 2.29e-07 \\  \hline
\multirow{2}*{E}
&1 &  14 & 5.35e-04 & 6.14e-04 & 2.90e-04 & 3.51e-04 \\  \cline{2-7}
&5 &  14 & 1.90e-05 & 1.04e-05 & 5.49e-06 & 4.79e-06 \\  \hline
\multirow{2}*{F}
&1 &  37 & 1.14e-05 & 4.49e-05 & 6.02e-06 & 2.16e-05 \\  \cline{2-7}  
&5 &  37 & 4.92e-07 & 8.19e-07 & 1.12e-07 & 2.60e-07 \\  \hline
 
\end{tabular}

\caption{LRA approximation of HSS matrices from \cite{XXG12}}
\end{center}
\label{tabhss}
\end{table}


\medskip


\noindent {\bf Acknowledgements:}
Our research was supported by NSF Grants CCF--1116736, CCF--1563942, and CCF--133834
and PSC CUNY Award 69813 00 48.




\begin{thebibliography}{\hspace{1cm}}


\bibitem[B00]{B00}
M. Bebendorf, Approximation of Boundary Element Matrices, {\em Numer. Math.},
{\bf 86,~4}, 565--589, 2000.


\bibitem[B10]{B10}
S. B\"{o}rm, {\em Efficient Numerical Methods
for Non-local Operators: $\mathcal H^2$-Matrix
Compression, Algorithms and Analysis},
European Math. Society, 2010.


\bibitem[B11]{B11}
M. Bebendorf,
Adaptive Cross Approximation of Multivariate Functions,
{\em Constructive approximation}, {\bf 34,~2}, 149--179, 2011.



\bibitem[BGH03]{BGH03}
S. B\"{o}rm, L. Grasedyck, W. Hackbusch, 
Introduction to Hierarchical Matrices with
Applications, 
{\em Engineering Analysis with Boundary Elements}, 
{\bf 27, ~5}, 405--422, 2003.
 

\bibitem[BGP05]{BGP05}
A. Bini, L. Gemignani, V. Y. Pan, 
Fast and Stable QR Eigenvalue Algorithms for Generalized Semiseparable Matrices and Secular Equation, 
{\em Numerische Mathematik}, {\bf 100,~3}, 373--408, 2005.




\bibitem[BR03]{BR03}
M. Bebendorf, S. Rjasanow,
Adaptive Low-Rank Approximation of Collocation Matrices,
{\em Computing}, {\bf 70,~1}, 1--24, 2003.


\bibitem[BY13]{BY13}
L. A. Barba, R. Yokota, 
How Will the Fast Multipole Method Fare in Exascale Era?
{\em SIAM News}, {\bf 46}, {\bf 6}, 1--3, July/August 2013.



\bibitem[C00]{C00}
B. A. Cipra, 
The Best of the 20th Century: Editors Name Top 10 Algorithms,
{\em SIAM News},  
{\bf 33,~4}, 2, 
May 16, 2000. 


\bibitem[CGR88]{CGR88}
J. Carrier, L. Greengard, V. Rokhlin,
A Fast Adaptive Algorithm for Particle Simulation,
{\em SIAM  Journal on Scientific Computing}, {\bf 9}, 669--686, 1988.


\bibitem[CGS07]{CGS07}
S. Chandrasekaran, M. Gu, X. Sun, J. Xia, J. Zhu,
A Superfast Algorithm for Toeplitz Systems of Linear Equations,
{\em SIAM. J. on Matrix Analysis and Applications}, {\bf 29}, {\bf 4}, 1247--1266, 2007.



\bibitem[CML15]{CML15}
 A. Cichocki, D. Mandic, L. D. Lathauwer, G. Zhou, Q. Zhao, C. Caiafa, H. A. Phan, Tensor Decompositions for Signal Processing Applications: From Two-Way to Multiway Component Analysis,
{\em IEEE Signal Processing Magazine}, {\bf 32,~2}, 145--163, March 2015.







\bibitem[EGH13]{EGH13}
Y. Eidelman, I. Gohberg, I. Haimovici,
{\em Separable Type Representations of Matrices and Fast Algorithms, volumes {\bf 1} and {\bf 2}},
Birkh{\"a}user, 2013.


\bibitem[GH03]{GH03}
 L. Grasedyck, W. Hackbusch, 
Construction and Arithmetics of H-Matrices,
{\em Computing}, {\bf 70,~4}, 295--334, 2003.


\bibitem[GOSTZ10]{GOSTZ10}
S. Goreinov, I. Oseledets, D. Savostyanov, E. Tyrtyshnikov, N. Zamarashkin,
How to Find a Good Submatrix, in 
{\em Matrix Methods: Theory, Algorithms, Applications} 
(dedicated to the Memory of Gene Golub, edited by V. Olshevsky and E. Tyrtyshnikov), 
pages 247--256,
World Scientific Publishing, New Jersey, ISBN-13 978-981-283-601-4, ISBN-10-981-283-601-2, 
 2010.


\bibitem[GR87]{GR87}
L. Greengard, V. Rokhlin,
A Fast Algorithm for Particle Simulation,
{\em Journal of Computational Physics}, {\bf 73}, 325--348, 1987.


\bibitem[GTZ97]{GTZ97}
 S. A. Goreinov, E. E. Tyrtyshnikov,  N. L. Zamarashkin, 
A Theory of Pseudo-skeleton Approximations,
{\em Linear Algebra and Its Applications}, {\bf 261}, 1--21, 1997.


\bibitem[HMT11]{HMT11}
N. Halko, P. G. Martinsson, J. A. Tropp,
Finding Structure with Randomness: Probabilistic Algorithms
for Constructing
 Approximate Matrix Decompositions, 
{\em SIAM Review}, {\bf 53,~2}, 217--288, 2011.



\bibitem[KS17]{KS17}
N. Kishore Kumar, J. Schneider,
Literature Survey on Low Rank Approximation of Matrices,
{\em Linear and Multilinear Algebra,} {\bf 65, 11}, 2212--2244, 2017, and
arXiv:1606.06511v1 [math.NA] 21 June 2016.



\bibitem[MRT05]{MRT05}
P. G. Martinsson, V. Rokhlin, M. Tygert, 
A Fast Algorithm for the Inversion of General
Toeplitz Matrices, 
{\em Comput. Math. Appl.}, {\bf 50}, 741--752, 2005.



\bibitem[O18]{O18}
A.I. Osinsky,
Rectangular Matrix Volume and Projective Volume Search Algorithms,
arXiv:1809.02334, September 17, 2018.



\bibitem[OZ18]{OZ18}
A.I. Osinsky, N. L. Zamarashkin,
 Pseudo-skeleton Approximations
 with Better Accuracy Estimates,
{\em Linear Algebra and Its Applications}, {\bf 537}, 221--249, 2018.


\bibitem[P15]{P15}
V. Y. Pan,
Transformations of Matrix Structures Work Again, 
{\em Linear Algebra and Its Applications}, 
{\bf 465}, 1--32, 2015.










\bibitem[T96]{T96}
E.E. Tyrtyshnikov,
 Mosaic-Skeleton Approximations,
{\em Calcolo}, {\bf 33,~1}, 47--57, 1996.


\bibitem[T00]{T00}
E. Tyrtyshnikov, Incomplete Cross-Approximation in the Mosaic-Skeleton Method,  
{\em Com\-put\-ing}, {\bf 64}, 367--380, 2000.


\bibitem[VVGM05]{VVGM05}
R. Vandebril, M. Van Barel, G. Golub, N. Mastronardi,
A Bibliography on Semiseparable Matrices,
{\em Calcolo}, {\bf 42,~3--4,} 249--270, 2005.


\bibitem[VVM07/08]{VVM07/08} 
R. Vandebril, M. Van Barel, N. Mastronardi,
{\em Matrix Computations and Semiseparable Matrices}
(Volumes 1 and 2),
The Johns Hopkins University Press, Baltimore, Maryland, 2007.


\bibitem[VVVF10]{VVVF10}
 M. Van Barel, R. Vandebril, P. Van Dooren, K. Frederix,
Implicit Double Shift $QR$-algorithm for Companion Matrices,
{\em Numerische Mathematik}, {\bf 116}, 177--212, 2010.


\bibitem[X12]{X12}
J. Xia, 
On the Complexity of Some Hierarchical Structured Matrix Algorithms,
{\em SIAM J. Matrix Anal. Appl.}, {\bf 33}, 388--410, 2012.


\bibitem[X13]{X13}
J. Xia, 
Randomized Sparse Direct Solvers,
{\em  SIAM J. Matrix Anal. Appl.}, {\bf 34}, 197--227, 2013.


\bibitem[XXCB14]{XXCB14}
J. Xia, Y. Xi, S. Cauley, V. Balakrishnan,
Superfast and Stable Structured Solvers for Toeplitz Least 
Squares via Randomized Sampling,
 {\em SIAM J. Matrix Anal. Appl.}, {\bf 35}, 44--72, 2014.


\bibitem[XXG12]{XXG12}
J. Xia, Y. Xi, M. Gu, 
A Superfast Structured Solver for Toeplitz Linear Systems via Randomized Sampling,
{\em SIAM J. Matrix Anal. Appl.}, {\bf 33}, 837--858, 2012.


\end{thebibliography}
\end{document}